\documentclass[a4paper,conference]{IEEEtran}
\IEEEoverridecommandlockouts
\usepackage[left=1.57cm,right=1.57cm,top=0.95cm,bottom=2.54cm]{geometry}
\usepackage{cite}
\usepackage[normalem]{ulem}
\usepackage{amsmath,amssymb,amsfonts}
\usepackage{algorithmic}
\usepackage{algorithm}
\usepackage{graphicx}
\usepackage{textcomp}
\usepackage{xcolor}

\newtheorem{lem}{Lemma}[section]
\newtheorem{theo}{Theorem}[section]
\newtheorem{coro}{Corollary}[section]

\newtheorem{rem}{Remark} 

\usepackage{dsfont}

\newtheorem{ass}{Assumptions}

\DeclareMathOperator*{\argmin}{arg\,min}
\newcommand{\proof}     {\paragraph*{Proof}}
\newcommand{\carre}     {\hfill$\Box$}

\def\BibTeX{{\rm B\kern-.05em{\sc i\kern-.025em b}\kern-.08em
    T\kern-.1667em\lower.7ex\hbox{E}\kern-.125emX}}
\begin{document}

\title{Online Optimization for Randomized Network Resource Allocation with  Long-Term Constraints
\thanks{This work was supported by a grant from the Natural Sciences and Engineering Research Council of Canada and Ericsson Canada. }
}

\author{\IEEEauthorblockN{ Ahmed Sid-Ali\IEEEauthorrefmark{1}, Ioannis Lambadaris\IEEEauthorrefmark{2}, Yiqiang Q. Zhao\IEEEauthorrefmark{1}, Gennady Shaikhet\IEEEauthorrefmark{1}, and Shima Kheradmand\IEEEauthorrefmark{2}}
\IEEEauthorblockA{\IEEEauthorrefmark{1}School of Mathematics and Statistics
\\}
\IEEEauthorblockA{\IEEEauthorrefmark{2}Department of Systems and Computer Engineering\\
Carleton University, Ottawa, Ontario\\
Emails: Ahmed.Sidali@carleton.ca; ioannis@sce.carleton.ca; zhao@math.carleton.ca; gennady@math.carleton.ca;\\ 
and  shimakheradmand@cunet.carleton.ca}

}

\maketitle

\begin{abstract}
This paper studies an online optimal resource reservation problem in communication networks. A typical network is composed of compute nodes linked by communication links. The system operates in discrete time and at each time slot, the administrator reserves server resources to meet future client requests. A cost is incurred for the made reservations. Then, after the reception of the job requests, some may be transferred from one server to another to be best accommodated at an additional transport cost. In particular, if certain job requests cannot be satisfied, a violation engenders a cost to pay for each blocked job. The goal is to minimize the overall reservation cost over finite horizons while maintaining the cumulative violation and transport costs under a certain budget limit. To tackle this problem, we propose a randomized procedure where the reservations are drawn randomly according to a sequence of probability distributions over the space of allowable reservations. The distributions are derived from an online saddle-point algorithm. For this proposed algorithm we analyze its performance by first establishing an upper bound for the associated \textit{$K$-benchmark regret}. We then develop an upper bound for the cumulative constraint violations. Finally, we present numerical experiments where we compare the performance of our algorithm with those of some basic deterministic resource allocation policies.  
\end{abstract}


\section{Introduction}

Online (convex) optimization (OCO) is a machine learning framework where a decision maker sequentially chooses a sequence of decision variables over time to minimize the sum of a sequence of (convex) loss functions. In OCO, the decision maker does not have full access to the data at once, but instead receives it incrementally over time. Thus, the decisions are made sequentially in an incomplete information environment. Moreover, in OCO, the data source is viewed as arbitrary, and thus no assumptions on the statistics of data sources
are made. This radically differs from classical statistical approaches such as Bayesian decision theory or Markov decision processes. Thus, only the observed data sequence's empirical properties matter, allowing us to address the dynamic variability of the traffic requests at modern communication networks. OCO has further applications in a wide range of fields, including control theory, game theory, economics, and signal processing, and has been extensively studied in the last decade; see, e.g. \cite{b12, b13, b14, b15} and the references therein for an overview. Moreover, given that the decision maker has only access to limited/partial information, globally optimal solutions are generally not realizable. Instead, one searches for algorithms that perform relatively well in comparison to the overall {\it ideal} best solution in hindsight which has full access to the data. This performance metric is referred to as \textit{regret} in the literature. In particular, if an algorithm incurs regret that increases sub-linearly with time, then one says that it achieves no regret. 

One important application of online optimization, as in the problem presented in our work, is in resource reservations that arise in situations where resources must be allocated in advance to meet future unknown demands. Typical examples arise in transportation systems, hotel reservations, or healthcare services handling systems. In this paper, we tackle a specific online resource reservation problem in a simple communication network topology. The network is composed of $N$ linked servers and the network administrator reserves resources at each server to meet future job requests. The specificity of this system is that one is allowed to transfer jobs from one server to another after receiving the job requests to better meet the demands. This couples the servers by creating dependencies between them. The reservation and transfer steps come with a cost that is proportional to the amount of resources involved. Moreover, if the job requests are not fully satisfied, a violation cost is incurred. The goal is then to minimize the reservation cost while maintaining the transfer and violation costs under a given budget threshold which leads to an online combinatorial optimization problem with long-term constraints.

To tackle the above-mentioned problem, we propose a randomized control policy that minimizes the cumulative reservation cost while maintaining the long-term average of the cumulative expected violation and transfer costs under the budget threshold. This transforms the online combinatorial optimization problem into an online {\it continuous} optimization problem on the space of probability distributions over the set of reservations. In particular, we propose an online saddle-point algorithm tailored to our model for which we derive an explicit upper bound for the incurred \textit{regret against $K$-Benchmark}, a concept introduced in \cite{b18} which allows us to examine the trade-offs between regret minimization and long-term budget constraint violations. 

The rest of the paper is organized as follows: in Section \ref{model-sec} we introduce the problem; in Section \ref{rand-formulation-sec}, we formalize the problem as a constrained online optimization problem on the simplex of probability distributions over the space of reservations; section \ref{algo-sec} contains our proposed online saddle point algorithm; then, in Section \ref{per-sec}, we present an upper bound for the $K$-benchmark regret together with the cumulative constraint violations upper bound, and finally, we present in Section \ref{num-sec} some numerical results where we compare the performance of our algorithm with some deterministic online reservation policies.  

\section{Online resource reservation in communication networks}  
\label{model-sec} 
Consider a network of $N$ servers connected by communication links. The system provides access to computing resources for clients. For simplicity, we suppose that there is a single type of resource (e.g. memory, CPU, etc.). Denote by $m_n$  the total number of resources available at the $n$-th server. Let $\mathcal{R}_n=\{1,\ldots,m_n\}$ be the set of possible reservations at the $n$-th server, and $\mathcal{R}=\prod_{n=1}^N\mathcal{R}_n$ be the set of possible reservations in the entire network. In particular, the system operates in discrete time slots $t=1,2,\ldots$ and, at each time slot $t$, the following processes take place chronologically: 
\begin{itemize} 
	\item \textbf{Resource reservation}: the network administrator selects the resources $A^t=(A^t_1,\ldots,A^t_N)\in\mathcal{R}$ to make available at each server.    
	\item \textbf{Job requests}: the network receives job requests $B^t=(B^t_1,\ldots,B^t_N)\in\mathcal{R}$ from its clients to each of its servers. 
	\item \textbf{Job transfer}: the network administrator can shift jobs between the servers to best accommodate the demand.  
\end{itemize} 
Let $\delta^t_{n,m}$ be the number of jobs transferred from server $n$ to server $m$ at time slot $t$ after receiving the job requests $B^t$. Notice that the coefficients $\delta^t_{n,m}$ depend on both $A^t$ and $B^t$. However, to keep the notation simple, we suppress this dependency in the sequel. 

Suppose that the reservations, the job transfers, and the violations incur costs defined as follows:

\paragraph*{\textbf{-Reservation cost}} given by the function
\begin{align*}
	C_R(A^t)=\sum_{n=1}^N  f^R_n(A^t_n).
\end{align*} 

\paragraph*{\textbf{-Violation cost}} incurred  when the job requests cannot be satisfied  
\begin{align*} 
	C_V(A^t,B^t)=\sum_{n=1}^N f^V_n\big(B_n^t-A_n^t-\sum_{m=1}^N\delta^t_{n,m}\big).
\end{align*} 

\paragraph*{\textbf{-Transfer cost}} incurred by the transfer of jobs between the different servers 
\begin{align*}
	&C_T(A^t,B^t)=\sum_{n=1}^N\sum_{m\neq n} f^T_n\big(\delta^t_{n,m}\big). 
\end{align*} 
Here, $f^R_n,f^V_n,f^T_n$, for $1\leq n\leq N$ are some given positive functions. Therefore, at each time slot $t$, and after receiving the job requests $B^t$, the job transfer coefficients $\{\delta^t_{n,m},1\leq n \leq N\}$ are solutions to the following offline minimization problem: 
\begin{equation}  
	\begin{split}
		\left\{
		\begin{tabular}{l}
			$\min \sum\limits_{n=1}^N\big(\sum\limits_{m\neq n}f^T_n(\delta^t_{n,m})+f^V_n\big(B_n^t-A_n^t-\sum\limits_{m\neq n}\delta^t_{n,m}\big)\big),$\\
			s.t. $\delta^t_{n,m}\leq\min\big\{ (B^t_n-A_n^t)^+,(A^t_{m}-B^t_{m})^+ \big\},\forall  n\neq m$. 
		\end{tabular}
		\right.
		\label{delta-opt-prob}
	\end{split}
\end{equation} 
The goal is then to minimize the reservation cost $C_R(A^t)$ while maintaining the sum  $C_V(A^t, B^t)+C_T(A^t, B^t)$ of the violation and transport costs under a given threshold $v>0$. Nevertheless, since at each slot $t$, the reservations $A^t$ are selected {\it before} the job request $B^t$ are received, this problem is of online optimization kind for which the optimal solution is out of reach. In particular, one cannot guarantee that the constraint is satisfied. Instead, one searches for online control policies that satisfy the constraints on average over large time horizons. More precisely, we consider the following online combinatorial constrained optimization problem with a long-term constraint:   
\begin{equation} 
	\begin{split}
		\left\{
		\begin{tabular}{l}
			$\min\limits_{\{A^t\}}\sum\limits_{t=1}^T C_R(A^t)$,\\
			s.t.:  $ \frac{1}{T} \sum\limits_{t=1}^T (C_V(A^t,B^t)+C_T(A^t,B^t))\leq v$,
		\end{tabular}
		\right.
		\label{opt-pb-on}
	\end{split}
\end{equation}
over finite time horizons $T\geq 1$. Notice that the sequence of job arrivals $\{B^t\}_{t\geq 1}$ is arbitrary, in the sense that we make no assumptions on its statistical properties which make the classical statistical inference methods unsuitable. Therefore, one searches for an online control policy that uses the cumulative information available so far to make the reservations at the next time slot. Of course, the goal is not to reach the optimal solution to $(\ref{opt-pb-on})$, and instead to obtain a total loss $\sum\limits_{t=1}^T C_R(A^t)$ that is not too large compared to some benchmark that knows the job requests in advance, and meanwhile, to ensure that the constraint is asymptotically satisfied on average over large time horizons. 
 
\begin{rem}
One can consider a different formulation of the problem by minimizing the total cost of reservation, violation, and transfer thus obtaining an unconstrained online optimization problem. This alternative formulation has been studied in \cite{b24} and a randomized exponentially weighted algorithm was proposed and analyzed. 
\end{rem} 
 
 \begin{figure}
 \center
 \includegraphics{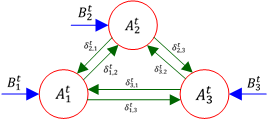}
 \caption{The network model for three nodes at time $t$}
 \end{figure}

\section{Online randomized reservations} 
\label{rand-formulation-sec}
Notice that the action set $\mathcal{R}=\prod_{n=1}^N\mathcal{R}_n$ is finite. Therefore, the classical online gradient type algorithms (see e.g. \cite{b17}) cannot be applied since no convexity assumption is possible. To overcome this issue, one possibility is to introduce randomization in the decision process. More precisely at every time step, the network administrator chooses a probability distribution $P^t$ over the set $\mathcal{R}$ of reservations, based on the past, and then draws his reservation $A^t$ at random according to $P^t$. In particular, denote by $\mathcal{P}(\mathcal{R})$ the space of probability distributions over the reservation set $\mathcal{R}$, and write $P=(p_a)_{a\in\mathcal{R}}$, then 
\begin{align*}   
	\mathcal{P}(\mathcal{R})=\bigg\{P\in\mathbb{R}^{|\mathcal{R}|}_{+}:\mbox{ $\sum\limits_{a\in\mathcal{R}}p_{a}=1$}\bigg\},
\end{align*} 
which is a convex, closed, bounded, and non-empty set. Define the expected reservation cost by 
\begin{align*}
	\mathbb{E}_{P^t}[C_R(A^t)]=\sum_{a\in\mathcal{R}}p^t_{a} C_R(a).
\end{align*}
Similarly, define the conditional expected transfer and violation costs given that the job request $B^t=b^t$, respectively, by
\begin{align*}
	\mathbb{E}_{P^t}[C_T(A^t,b^t)]=\sum_{a\in\mathcal{R}}p^t_{a}C_T(a,b^t)
\end{align*}
and 
\begin{align*}
	\mathbb{E}_{P^t}[C_V(A^t,b^t)]=\sum_{a\in\mathcal{R}}p^t_{a}C_V(a,b^t).
\end{align*}  
The goal is then to find a randomized control policy that behaves relatively well in comparison to the best solutions in hindsight. In particular, let $A_K^*$ be the optimal static solution to the constrained combinatorial optimization problem  
\begin{equation} 
	\begin{split}
		\left\{
		\begin{tabular}{l}
			$\min\limits_{A\in\mathcal{R}}  C_R(A)$,\\
			s.t.: $\sum\limits_{k=t}^{t+K-1}(C_V(A,b^{k})+C_T(A,b^{k}))\leq Kv,$\\
			$\qquad\qquad\qquad\qquad\qquad\forall t=1,\ldots,T-K+1$,
		\end{tabular}
		\right.
	\end{split}
	\label{opt-pb-hind-no-rand}
\end{equation}
if the values of job requests $\{B^t\}_{t=1}^T=\{b^t\}_{t=1}^T$ were known in advance. Moreover, let $\{A_t\}_{t=1}^T$ be the sequence of random reservations generated according to the sequence $\{P^t\}_{t=1}^T$ of probability distributions. Therefore, the regret of not playing $A^*_K$ over the horizon $T$ is given by
\begin{align}
	\tilde{R}^{K}_T= \sum_{t=1}^T \big( C_R(A^t)- C_R(A_K^*)\big). 
	\label{reg-non-rand}
\end{align}

The definition in $(\ref{reg-non-rand})$ is known in the literature as the ``regret against $K$-benchmark"; see, e.g. \cite{b18,b20}. The idea is that by varying $K$, this metric provides sufficient flexibility to study the question of adversarial long-term constraints. It is well known from \cite{b21} that for $K=T$, sub-linear regret is not attainable in a general adversarial setting. In counterpart, for $K=1$, a sub-linear regret is achievable under some reasonable assumptions as we will see in the sequel. Therefore, we present in this paper a general upper bound for the $K$-benchmark regret which is sub-linear in $T$ when $K=o(\sqrt{T})$. Furthermore, the regret $\tilde{R}^{K}_T$ in $(\ref{reg-non-rand})$ is a random variable since the reservations $\{A_t\}_{t=1}^T$  are generated randomly. Alternatively, one can introduce the deterministic regret $R^K_T$ defined in terms of expected costs. More precisely, let the probability distribution $P^K_*\in\mathcal{P}(\mathcal{R})$ be the optimal solution to the following optimization problem:    
	\begin{equation}
		\begin{split}
			\left\{
			\begin{tabular}{l}
				$\min\limits_{P\in\mathcal{P}(\mathcal{R})}\sum\limits_{t=1}^T \mathbb{E}_{P}[C_R(A^t)]$,\\
				s.t.: $\sum\limits_{k=t}^{t+K-1}\mathbb{E}_{P}[C_V(A^{k},b^{k})+C_T(A^{k},b^{k})]\leq Kv,$\\
				$\qquad\qquad\qquad\qquad\qquad\forall t=1,\ldots,T-K+1$,
			\end{tabular}
			\right.
			\label{p-star-stat-K}
		\end{split}
	\end{equation}   
	for finite $1\leq K\leq T$, if the values of job requests $\{B^t\}_{t=1}^T=\{b^t\}_{t=1}^T$ were known in advance. Then, one defines the deterministic \textit{regret} of not playing according to $P^K_*$ over the horizon $T$ as 
	\begin{equation}
		\begin{split}
			&R^K_T= \sum\limits_{t=1}^T \big( \mathbb{E}_{P^t}[C_R(A^t)] - \mathbb{E}_{P_*^K}[C_R(A^t)]\big). 
			\label{stat-reg-K}
		\end{split}
	\end{equation} 
	In Section \ref{per-sec}, we establish upper-bounds for $R^K_T$ and $\tilde{R}^{K}_T$. In addition to the regret, we analyze the expected cumulative constraint violations over finite horizons $T>0$,
\begin{align} 
	\Upsilon_T= \sum_{t=1}^T\big( \mathbb{E}_{P^t}[C_V(A^t,b^t)+C_T(A^t,b^t)]-v \big).
	\label{cumul-const}
\end{align} 

Let us suppose the following assumptions:  

\begin{ass}
	\label{ass}
\begin{enumerate}
	\item Suppose that there exists a constant $\Theta>0$ such that, for all $a,b\in\mathcal{R}$,
	\begin{align*}
		\big|C_R(a)\big|\leq \Theta,\big|C_T(a,b)\big|\leq \Theta,\mbox{ and } \big|C_V(a,b)\big|\leq \Theta. 
	\end{align*}
	\item There exists a subspace $\mathcal{E}\subset\mathcal{P}(\mathcal{R})$ such that, for all $P\in\mathcal{E}$,  
	\begin{align*}  
		\mathbb{E}_P[C_V(A,b)+C_T(A,b)]-v\leq 0, \mbox{ for all $b\in\mathcal{R}$}. 
	\end{align*}
	\item There exist a probability distribution $\tilde{P}\in \mathcal{E}$, called a \textit{Slater} distribution, and a constant $\eta>0$ such that,   
	\begin{align*}
		\mathbb{E}_{\tilde{P}}[C_V(A,b)+C_T(A,b)]-v\leq -\eta, \mbox{ for all $b\in\mathcal{R}$}. 
	\end{align*}
\end{enumerate} 
\end{ass}

\begin{rem}
	The subspace $\mathcal{E}$ represents the feasible set in which the constraints are satisfied at all slots $t\geq 1$. The Slater vector $\tilde{P}$ is an interior point for $\mathcal{E}$ and its existence is a classical sufficient condition for strong duality to hold for a convex optimization problem; see, e.g. \cite{b12}[Prop. 3.3.9]. 
\end{rem}  
  
\section{Online Lagrange multipliers approach}   
\label{algo-sec}
A classical approach for solving constrained convex optimization problems is the Lagrange multipliers method. This relies in particular on the Karush-Kuhn-Tucker theorem which states that finding an optimal point for the constrained optimization problem is akin to finding a saddle point of the Lagrangian function (\textit{cf.} \cite[Section I.5]{b3}). Popular saddle-point methods are the Arrow-Hurwicz-Uzawa algorithms introduced in \cite{b2} and widely used in the literature; see, e.g. \cite{b5,b6,b7} for an overview. These methods alternate a minimization of the Lagrangian function with respect to the primal variable and a gradient ascent with respect to the dual variable given the primal variable. Online versions of these methods have also been proposed recently (see, e.g., \cite{b8,b9,b16}). In the same spirit as the aforementioned references, we develop here an online saddle-point algorithm to solve the problem detailed in Section \ref{rand-formulation-sec}. To this end, let us first introduce the sequence of per slot optimization problem 
\begin{equation}
\begin{split}
\left\{
\begin{tabular}{l}
$\min\limits_{P^t\in\mathcal{P}(\mathcal{R})} \mathbb{E}_{P^t}[C_R(A^t)]$,\\
s.t.: $\mathbb{E}_{P^t}[C_V(A^t,b^t)+C_T(A^t,b^t)]-v\leq 0$,
\end{tabular}
\right.
\end{split}
\end{equation}
and the corresponding instantaneous online Lagrangian functions
\begin{equation}
\begin{split}
\mathcal{L}^t(P^t,\lambda_t)&=\mathbb{E}_{P^t}[C_R(A^t)] \\
&\qquad+\lambda_t \big(\mathbb{E}_{P^t}\big[C_V(A^t,b^t)+C_T(A^t,b^t)\big]-v \big),
\end{split}
\end{equation} 
where $\{\lambda_t\}_{t\geq 1}$ is a sequence of Lagrange multipliers and $P^t\in\mathcal{P}(\mathcal{R})$ for all $t\geq 1$.
 Now, since at each time slot $t$, the probability distribution $P^t$ is selected before the job request $b^t$ is observed, $\mathcal{L}^t(P^t,\lambda_t)$ cannot be computed beforehand. Nonetheless, one uses the most recent information $B^{t-1}=b^{t-1}$ available to build the following sequence of approximate Lagrangian functions:  
 \begin{equation}
 \begin{split}
\tilde{\mathcal{L}}^t(P^t,\lambda_t)&= \mathbb{E}_{P^t}[C_R(A^t)]\\
&\quad+\lambda_t \big(\mathbb{E}_{P^t}[C_V(A^t,b^{t-1})+C_T(A^t,b^{t-1})\big]-v \big).
 \end{split}
 \label{appr-lag}
 \end{equation} 
By noticing that the function $\mathcal{L}^t$ is linear in $P^t$, it is easy to see that our approximation is equivalent to the ones proposed in \cite{b8} and \cite{b17}. Note that the authors in \cite{b17} weight the objective function by an additional penalty parameter. 
 
Observe that the minimization of $(\ref{appr-lag})$ at each slot $t$ could result in an oscillatory behavior for $P^t$, especially if the job requests are not smooth. To overcome this, one introduces an additional proximal term and then, at each time slot $t$, selects $P^t$ that solves the following minimization problem: 
\begin{equation}
\begin{split}
\min_{P\in\mathcal{P}(\mathcal{R})}\bigg(\tilde{\mathcal{L}}^t(P,\lambda_t)+\frac{1}{2\alpha} \big\|P-P^{t-1}\big\|^2\bigg),
\label{appro-lag-quad}
\end{split}
\end{equation} 
with $\alpha>0$ is a positive scalar and $\|\cdot\|$ stands for the Euclidean distance, thus obtaining a proximal point approach. In particular, the proximal methods are well-known in the classical optimization literature where their main advantage is to transform the original objective function into a strongly convex one. Therefore, the convergence of these methods does not require strict convexity; see, e.g. \cite{b4,b10,b11} and the references therein for more details. For online optimization algorithms, the proximal term $ \big\|P-P^{t-1}\big\|^2$ acts as a learning term that accumulates knowledge from the past allowing to improve the performance of the algorithms in terms of regret. Thus many authors used a similar term in their algorithms; see, e.g. \cite{b8,b17,b19}. Finally, the Lagrange multiplier $\lambda_t$ is in turn updated at each step $t$ using a projected gradient ascent step. Namely, 
\begin{align} 
\lambda_{t+1}&=\bigg(\lambda_t+\mu \big( \mathbb{E}_{P^t}\big[C_V(A^t,b^{t-1})+ C_T(A^t,b^{t-1})\big]-v \big)\bigg)^+,
\label{lambd-update}
\end{align}
where $\mu>0$ is the step-size, and, for any $x\in\mathbb{R}$, $(x)^+=\max\{0,x\}$. Notice that the update in $(\ref{lambd-update})$ is based on the gradient of the approximate Lagrangian function $\tilde{\mathcal{L}}^t$ with respect to $\lambda_t$ instead of the gradient of the Lagrangian function $\mathcal{L}^t$ as proposed in \cite{b8} and \cite{b19}. More precisely, we update $\lambda_{t+1}$ using $b^{t-1}$ even though $b^t$ is available after the reservation is made. The reason for that is to use the same function in the primal and dual update equations which allows us to obtain an upper bound for the regret independent of the sequence of job arrivals. Moreover, by the linearity of $\tilde{\mathcal{L}}^t$ with respect to $P^t$, the update in $(\ref{lambd-update})$ amounts to the virtual queue update proposed in \cite{b17,b18} with the difference that we include a step-size parameter $\mu$.  

We summarize our approach in Algorithm \ref{sad-algo}.
\begin{algorithm}
\caption{Randomized Saddle-Point Algorithm for Online Resource Reservation}
\begin{itemize}
\item Initialize the values of $P^0\in\mathcal{P}(\mathcal{R})$, $b^0\in\mathcal{R}$, and $\lambda_1\geq 0$ 
\item At each step $t=1,2,\ldots$  
\begin{itemize}
\item Observe $b^{t-1}$ and chose $P^t$ such that:
\end{itemize}
\begin{align*}
P^t=\argmin_{P\in\mathcal{P}(\mathcal{R})}\big( \tilde{\mathcal{L}}^t(P,\lambda_t)+\frac{1}{2\alpha} \|P-P^{t-1}\|^2\big)
\end{align*} 
\begin{itemize}
\item Update $\lambda_t$ by:
\end{itemize}
\end{itemize}   
\begin{align*}
\lambda_{t+1}&=\big(\lambda_t+\mu \big( \mathbb{E}_{P^t}\big[C_V(A^t,b^{t-1})+ C_T(A^t,b^{t-1})\big]-v \big)\big)^+
\end{align*} 
\label{sad-algo}
\end{algorithm}

\section{Performance analysis}
\label{per-sec}
 We analyze here the performance of Algorithm \ref{sad-algo} over finite time horizons $T>0$, first in terms of cumulative constraint violations and then in terms of regret. To facilitate the reading, all proofs are postponed to the appendices. 

\subsection{Violation constraint bound}
We start by establishing an upper bound for the cumulative constraint violations $\Upsilon_T$ defined in $(\ref{cumul-const})$. 
\begin{theo} 
Let $\lambda_1=0$. Then, for any positive integer $\aleph\in\mathbb{N}$, the Lagrange multiplier $\lambda_t$ given by the update equation in $(\ref{lambd-update})$ is bounded as: $
\lambda_t\leq \bar{\lambda}:= \theta \aleph,\mbox{ for all $t\in\{1,2,\ldots\}$}$,
where $\theta=\max\big\{\varrho,\chi \big\}$, with
\begin{align*}
\chi&= \frac{1}{\eta\aleph}\bigg(  \frac{\mu \big(4\Theta^2+v^2 \big)}{2} +  \frac{\alpha \Theta^2}{4}+\frac{1}{2\alpha(\aleph+1)} +  \Theta\bigg)\\
&\qquad+ \frac{\varrho(\aleph+2)}{2\aleph}, 
\end{align*}
and $\varrho=\mu (2\Theta-v)$, $\alpha,\mu$ are the step sizes, $\Theta,\eta$ are defined as previously, and $v$ is the violation and transport cost threshold. Moreover, for any $T>0$, the cumulative violations related to Algorithm \ref{sad-algo} satisfies:  $\Upsilon_T \leq \frac{\theta \aleph}{\mu}$.
\label{Theo-fit} 
\end{theo}  
\proof See Appendix \ref{proof theo 1}.

\subsection{Regret against $K$ benchmark}
We first establish an upper bound for the deterministic regret $R^K_T$ defined in $(\ref{stat-reg-K})$ associated with the sequence $\{P^t\}_{t=1}^T$ of probability distributions drawn from Algorithm \ref{sad-algo}.
  \begin{theo}
   Let $\lambda_1=0$. Then, the deterministic regret associated with the sequence of probability distributions generated by Algorithm \ref{sad-algo} satisfies, for any $T>0$ and $1\leq K\leq T$,
\begin{align*}
R^K_T&\leq \mu(2\Theta+v)^2 \frac{(K-1)(2K-1)}{6} \\
 &\quad+(T-K)\bigg(\frac{1}{2} K\big(4 \Theta^2+v^2 \big)\mu + \frac{ \Theta^2}{4} \alpha\bigg) \\
 &\quad+  K(K-1)\Theta+\frac{1}{\alpha}.
 \label{reg-up-bound-K} 
\end{align*}
\label{Theo-regret-K}
\end{theo} 
\proof See Appendix \ref{K-bench-proof}

The next corollary establishes an upper bound for the regret associated with the corresponding sequence $\{A^t\}_{t=1}^T$ of randomized reservations.
\begin{coro}
 Let $\lambda_1=0$. Then, for any $0<\delta <1$, the regret $\tilde{R}^K_T$ related to Algorithm \ref{sad-algo} satisfies, for any $T>0$ and $1\leq K\leq T$,
 		\begin{equation}
 	\begin{split}
 	\tilde{R}^K_T &\leq R^{K}_T+\sqrt{2 \log (\delta^{-1}) T \Theta^2}
 	\end{split}
 	\label{reg-non-rand-bound}
 \end{equation}
 with probability at least $1-\delta$. 
 \label{coro}
	\end{coro}
\proof See Appendix \ref{coro-proof}.

\begin{rem}
Notwithstanding that the regret in Theorem \ref{Theo-regret-K}, and thus in Corollary \ref{coro}, is not sub-linear in $T$, one can choose the parameters $\alpha$ and $\mu$ such that for $T$ large enough, the regret becomes very small when $K=o(\sqrt{T})$. Indeed, fix for instance $K=1$ and $\varepsilon>0$. Moreover, let us chose $\alpha=\varepsilon^{\iota},\mu=\varepsilon^{\gamma}$ and $T'=\lceil\big( 1/\varepsilon\big)^{\kappa} \rceil$. Therefore, from Theorem \ref{Theo-regret-K} one obtains, for any $T>T'$,  
   \begin{equation} 
   \begin{split}  
    \frac{R^K_T}{T} &\leq \mathcal{O}(\varepsilon^{\gamma}+\varepsilon^{\iota}+\varepsilon^{\kappa-\iota}),
   \end{split}  
   \end{equation}
   and similarly from Theorem \ref{Theo-fit}
   \begin{equation}
   \frac{\Upsilon_T}{T} \leq   \varepsilon^{\iota+\kappa-\gamma}+ \varepsilon^{\kappa-\iota-\gamma}.  
   \end{equation}
   In particular, if $\iota=\gamma=1$ and $\kappa=3$, then 
   \begin{equation} 
   \begin{split}  
    \frac{R^K_T}{T} &\leq \mathcal{O}(\varepsilon+\varepsilon^{2}), \mbox{ and } \frac{\Upsilon_T}{T} \leq  \mathcal{O}( \varepsilon+ \varepsilon^{3}). 
   \end{split}  
   \end{equation}
   This leads us to the concept of convergence time to an $\mathcal{O}(\varepsilon)$-approximation introduced in \cite{b17}.
\end{rem}

\section{Numerical experiments}
\label{num-sec}
We propose in this section to test our proposed approach on some synthetic data. Moreover, we compare its performance against the following simple deterministic  policies:  
 
\paragraph{Lazy bang-bang policy} this policy is described in Algorithm \ref{lazy-algo} 
\begin{algorithm}
\caption{Lazy bang-bang policy}
\begin{itemize}
\item For each $t=1,\ldots,T$, observe $b^{t-1}$ and let $A^t=b^{t-1}$
\end{itemize}
\label{lazy-algo}
\end{algorithm}

\paragraph{Naive bang-bang policy}  This policy is detailed in Algorithm $\ref{naive-algo}$. 

\begin{algorithm}
\caption{Naive bang-bang policy}
\begin{itemize} 
\item At each slot $t=1,\ldots,T$ observe $b^{t-1}$ and select $A^t$ that solves: 
\begin{align*}
\left\{
\begin{tabular}{l}
$\min\limits_{A\in\mathcal{A}}C_R(A)$\\
$\mbox{subject to: } C_V(A,b^{t-1})+C_T(A,b^{t-1}) \leq v$
\end{tabular}
\right.
\end{align*}
\end{itemize}
\label{naive-algo}
\end{algorithm}

\paragraph{Lagrangian deterministic algorithm} this policy, detailed in Algorithm \ref{queu-algo}, mimics our saddle-point algorithm with the difference that the minimization is done over the set of reservations instead of  the simplex of probability distributions.   
\begin{algorithm}[h]
\caption{Lagrangian combinatorial algorithm}
\begin{itemize}
\item Set $\lambda_1=0$
\item At each slot $t=1,\ldots,T$, observe $b^{t-1}$ and select $A^t$ as:
\end{itemize}
\begin{align*}
&\argmin\limits_{A\in\mathcal{A}}\big[C_R(A)+\lambda_t\big(C_V(A,b^{t-1})+C_T(A,b^{t-1})-v \big)  \big]\\
&\lambda_{t+1}=\big(\lambda_t+ \big(C_V(A^t,b^{t-1})+C_T(A^t,b^{t-1})-v \big) \big)^+
\end{align*}
\label{queu-algo}
\end{algorithm}

\begin{figure}[ht]
\centering
\includegraphics[scale=0.38]{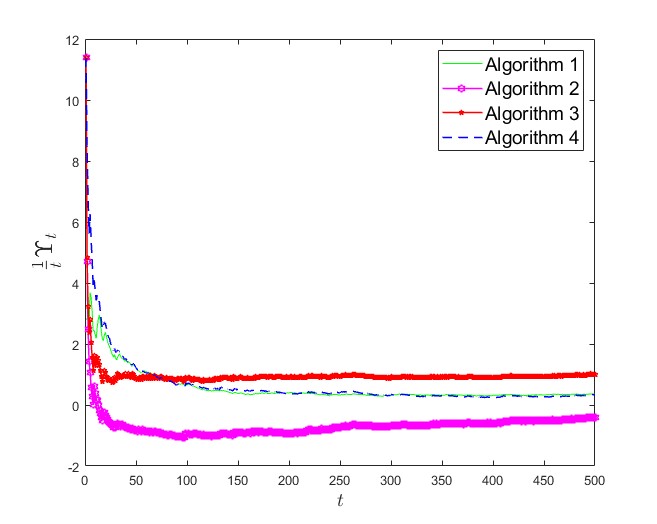}
\caption{Time average constraint violations}
\label{viol-fig}
\end{figure}

\begin{figure}[ht]
\centering
\includegraphics[scale=0.38]{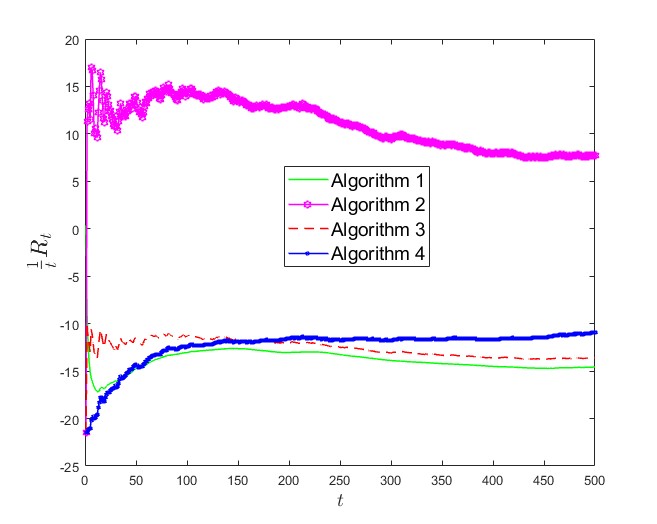}
\caption{Time average $1$-Benchmark regret}
\label{reg-1-fig}
\end{figure}

\begin{figure}[ht]
\centering
\includegraphics[scale=0.38]{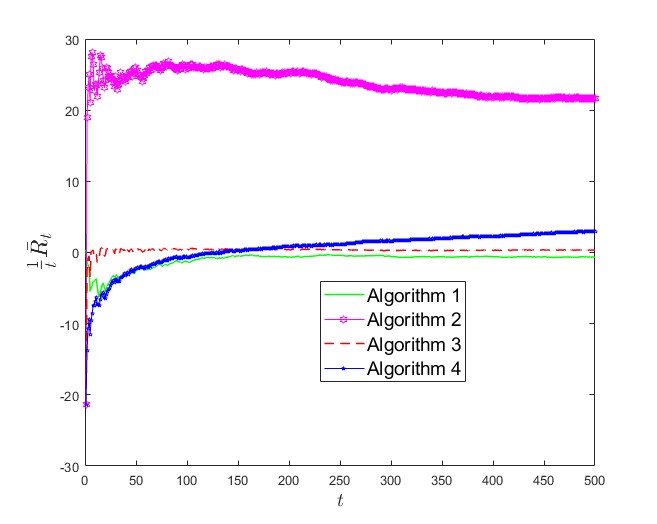}
\caption{Time average $T$-Benchmark regret}
\label{reg-T-fig}
\end{figure}

We then test the four algorithms on synthetic data generated on a simple two-server network linked by a bi-directional link. The algorithms are run for a total of $T=500$ time slots using the following set of parameters: $n_1=7,n_2=8,v=2$, and  $\mu=0.1$. In addition, we fix the value of the penalty parameter at $\alpha=0.001$ in Figures \ref{viol-fig},  \ref{reg-1-fig} and \ref{reg-T-fig}. In Figure \ref{conv-fig}, we run Algorithm \ref{sad-algo} with two different values $\alpha=0.01$ and $\alpha=0.001$ to test the effect of this parameter on the Euclidean distance $\|P^t-P{t-1}\|$ between consecutive terms in the sequence $\{P^t\}_{t=1,\ldots,500}$ of probability distributions generated by running Algorithm \ref{sad-algo}. Also, the cost functions are given, for all $x>0$ by: $f^1_R(x)=(0.3)x^2$, $f^2_R(x)=(0.1)x^3$, $f^1_V(x)=0.1 x^2$,  $f^2_V(x)=0.2 x^2$, $f^1_T(x)=log(x+1)$, and  $f^2_T(x)=log((x+1)/2$ and are equal to $0$ for $x<0$. 

From Figure \ref{viol-fig}, one can observe that Algorithm $\ref{sad-algo}$ and Algorithm $\ref{queu-algo}$ tend to converge in terms of time-average cumulative constraint violations, this is due to the minimization of the Lagrangian functions that take into account, at each time slot, the cumulative constraint violations so far and thus learn from the past. On the other hand,  Algorithm \ref{lazy-algo} and Algorithm \ref{naive-algo}, which use, at each time slot, solely the information from the previous slot, do not show tendency and thus follow the fluctuations of the job arrivals. In Figure \ref{reg-1-fig}, we plot the time average $1$-Benchmark regret. Here we can see that Algorithm \ref{sad-algo} outperforms the other algorithms. Furthermore, in Figure \ref{reg-T-fig} we plot the time-average $T$-benchmark regret. Surprisingly, even though the theoretical bound, in this case, is quadratic in $T$, the result shows a sub-linear regret for Algorithm \ref{sad-algo}. This opens the fundamental question of under which condition on the sequence of job requests one can hope to obtain a theoretical sub-linear upper bound for the $T$-benchmark regret. This, however, goes beyond the scope of this paper and will be the subject of future investigations. Finally, we plot in Figure \ref{conv-fig} the Euclidean distance between the consecutive probability distributions generated by Algorithm \ref{sad-algo} with two different values of the penalty parameter $\alpha$. Smaller values of $\alpha$ produce fewer fluctuations in the sequence $\{P^t\}_{t\geq 1}$ since it puts more weight on the quadratic Tikhnov term in the minimization.    

\begin{figure}[ht]
\label{conv-fig}
\centering
\includegraphics[scale=0.38]{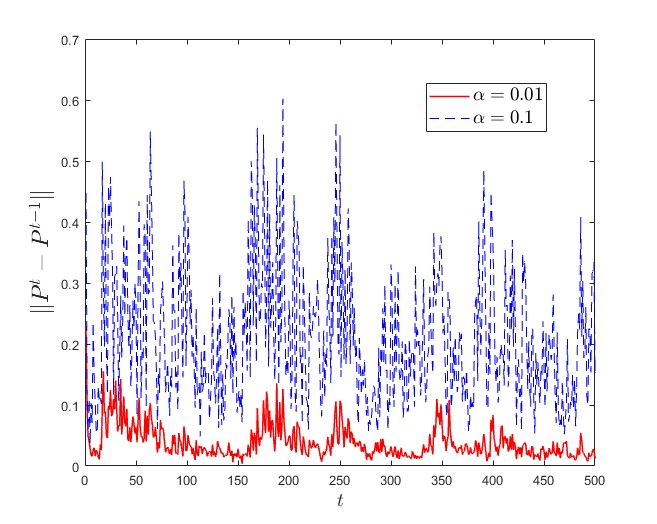}
\caption{Effect of $\alpha$ parameter on the distance between successive probability distributions}
\label{conv-fig}
\end{figure}

\section{Conclusion}
In this work, we studied an online optimal resource reservation problem in communication networks. The coupling between the servers, introduced by the option to transport jobs between servers at a cost, creates a more complex problem since the job requests arriving at each server together with the reservation made on it affect the decision on the other servers. Adding randomization allows us to embed the problem into the continuous space of probability distributions and thus apply a saddle point algorithm for which we have presented  $K$-benchmark regret and cumulative violations upper bounds. Moreover, the numerical experiments showed the comparison between our algorithm and some simple deterministic policies. The study leads to several interesting questions. The established $K$-benchmark regret upper bound is sublinear in time for $K=o(\sqrt{T})$, however, our preliminary simulations showed that the regret seems to be sublinear in time even for $T$-Benchmark. Thus, one may wonder under which conditions on the job requests sequence one can obtain sublinearity for $T$-benchmark regret. In particular, it is well known that in a general adversarial environment, sublinearity, in this case, is not possible. However, in many networking applications, the job request sequence may not be adversarial.

\appendix 
\subsection{Technical results}
Define the one-slot drift as $\Delta(\lambda_t)=\frac{1}{2}\big(\lambda^2_{t+1}-\lambda_t^2\big)$. Using $(\ref{lambd-update})$ one deduces the following lemma: 
\begin{lem}
	At each slot $t=1,2,\ldots$, one has:  
	\begin{align*}  
		\Delta(\lambda_t) &\leq  B + \mu  \lambda_t \big(\mathbb{E}_{P^t}[C_V(A^t,b^{t-1})+C_T(A^t,b^{t-1})]-v\big),
	\end{align*}
	where $B=\frac{1}{2} \mu^2 \big(4\Theta^2+v^2 \big)$.
	\label{lag-mult-bound}
\end{lem}

\begin{lem}
	For every vector $P \in \mathcal{P}(\mathcal{R})$ and every slot $t \in\{1, 2,\ldots\}$, the solution $P^t$ given by Algorithm $\ref{sad-algo}$
	satisfies 
	\begin{equation}
		\begin{split}  
		\frac{\Delta(\lambda_t)}{\mu}&\leq \frac{B}{\mu}+ \mathbb{E}_{P}[C_R(A^t)]-\mathbb{E}_{P^{t-1}}[C_R(A^{t-1})]\\
		&+\lambda_t \big(\mathbb{E}_{P}[C(A^t,b^{t-1})\big]-v \big)\\
		&+\frac{1}{2\alpha} \bigg(\|P-P^{t-1}\|^2-\|P-P^{t}\|^2\bigg)+\frac{\alpha \Theta^2}{4}. 
		\label{delta-bound-ineq}
	\end{split}
	\end{equation}
	\label{Delta-bound}
\end{lem}
\proof Let us define, for all $t$, the sequence of functions  
\begin{align*}  
	h_t(P)=\tilde{\mathcal{L}}^{t}(P,\lambda_{t})+\frac{1}{2\alpha} \|P-P^{t-1}\|^2. 
\end{align*}
It is easy to see that $h_t$ is $\frac{1}{\alpha}$-strongly convex. Moreover, recall that the probability distribution $P^{t}$ given by Algorithm \ref{sad-algo} at time slot $t$ minimizes $h_t(P)$.  Therefore, for any $P\in\mathcal{P}(\mathcal{R})$,
\begin{align*}
	h_t(P^{t})\leq h_t(P)-\frac{1}{2\alpha}\|P-P^{t}\|^2.
\end{align*}  
Therefore,
\begin{align*}  
	&\tilde{\mathcal{L}}^{t}(P^t,\lambda_{t})+\frac{1}{2\alpha} \|P^t-P^{t-1}\|^2\\
	&\leq \tilde{\mathcal{L}}^{t}(P,\lambda_{t})+\frac{1}{2\alpha} \|P-P^{t-1}\|^2-\frac{1}{2\alpha}\|P-P^{t}\|^2.
\end{align*}
Recall that 
\begin{equation}
	\begin{split}
		\tilde{\mathcal{L}}^t(P^t,\lambda_t)&= \mathbb{E}_{P^t}[C_R(A^t)]+\lambda_t \big(\mathbb{E}_{P^t}[C(A^t,b^{t-1})\big]-v \big),
	\end{split}
\end{equation}  
where $C(a,b):=C_V(a,b)+C_T(a,b)$. Then 
\begin{align*}  
	&\mathbb{E}_{P^t}[C_R(A^t)]+\lambda_t \big(\mathbb{E}_{P^t}[C(A^t,b^{t-1})\big]-v \big)\\
	&\qquad+\frac{1}{2\alpha} \|P^t-P^{t-1}\|^2\\
	&\leq \mathbb{E}_{P}[C_R(A^t)]+\lambda_t \big(\mathbb{E}_{P}[C(A^t,b^{t-1})\big]-v \big)\\
	&\qquad+\frac{1}{2\alpha} \|P-P^{t-1}\|^2-\frac{1}{2\alpha}\|P-P^{t}\|^2. 
\end{align*}  
By Lemma \ref{lag-mult-bound} one gets 
\begin{align*}  
	&\mathbb{E}_{P^t}[C_R(A^t)]+ \frac{\Delta(\lambda_t)}{\mu}+\frac{1}{2\alpha} \|P^t-P^{t-1}\|^2\\
	&\leq \frac{B}{\mu}+ \mathbb{E}_{P}[C_R(A^t)]+\lambda_t \big(\mathbb{E}_{P}[C(A^t,b^{t-1})\big]-v \big)\\
	&+\frac{1}{2\alpha} \|P-P^{t-1}\|^2-\frac{1}{2\alpha}\|P-P^{t}\|^2.
\end{align*}   
Adding and subtracting $\mathbb{E}_{P^{t-1}}[C_R(A^{t-1})]$ to the left hand side of the last inequality gives   
\begin{align*}  
	&\mathbb{E}_{P^{t-1}}[C_R(A^{t-1})]+\mathbb{E}_{P^t-P^{t-1}}[C_R(A^t)]\\
	&\qquad+ \frac{\Delta(\lambda_t)}{\mu}+\frac{1}{2\alpha} \|P^t-P^{t-1}\|^2\\
	&\leq \frac{B}{\mu}+ \mathbb{E}_{P}[C_R(A^t)]+\lambda_t \big(\mathbb{E}_{P}[C(A^t,b^{t-1})\big]-v \big)\\
	&\qquad+\frac{1}{2\alpha} \|P-P^{t-1}\|^2-\frac{1}{2\alpha}\|P-P^{t}\|^2.
\end{align*} 
Simple calculations using the assumptions above give 
\begin{align*}
	-  \mathbb{E}_{P^{t}-P^{t-1}}[C_R(A^t)]-\frac{1}{2\alpha} \|P^{t}-P^{t-1}\|^2\leq \frac{\alpha \Theta^2}{4}, 
\end{align*}
from which one deduces $(\ref{delta-bound-ineq})$.
\carre

\subsection{Proof of Theorem \ref{Theo-fit}}
\label{proof theo 1}
We first establish an upper bound for the Lagrange multipliers $\lambda_t$.
\begin{lem}
	Suppose that Slater's condition holds. Therefore, for any positive integer $\aleph\in\mathbb{N}$, the Lagrange multiplier $\lambda_t$ given by the update equation $(\ref{lambd-update})$ is bounded as follows:
	\begin{align*} 
		\lambda_t\leq \theta \aleph,\mbox{ for all $t\in\{1,2,\ldots\}$},
	\end{align*}  
	where 
	\begin{align}
		\theta=\max\{\varrho,\digamma \}
		\label{thet-def}
	\end{align}
	with  $\varrho=\mu (2\Theta-v)$ and 
	\begin{align*}
\digamma=\frac{1}{\eta\aleph}\bigg(  \frac{B}{\mu}+  \frac{\alpha \Theta^2}{4}+\frac{1}{(\aleph+1)}  \frac{1}{2\alpha} +  \Theta\bigg)+ \frac{\varrho(\aleph+2)}{2\aleph}.
	\end{align*} 
	\label{queu-bound-lem}
\end{lem}

\proof suppose that $Q^*$ is Slater's vector. Therefore, by the inequality in $(\ref{delta-bound-ineq})$ one obtains 
\begin{equation} 
	\begin{split}    
		\frac{\Delta(\lambda_t)}{\mu}&\leq \frac{B}{\mu}+ \Theta-\eta \lambda_t+\frac{\alpha \Theta^2}{4}\\
		&+\frac{1}{2\alpha} \bigg(\|Q^*-P^{t-1}\|^2-\|Q^*-P^{t}\|^2\bigg). 
		\label{Delt-bound}
	\end{split}
\end{equation} 
Moreover, from the queue update equation in $(\ref{lambd-update})$ and the assumptions above one gets that, for all $t\in\{1,2,\ldots\}$, 
\begin{align}
	\lambda_{t+1}-\lambda_t \leq \mu (2\Theta-v). 
	\label{lamb-succ-bound}
\end{align} 
Therefore, for all $t\in\{1,2,\ldots\}$,
\begin{align*}
	\lambda_{t} \leq (\mu (2\Theta-v))t. 
\end{align*} 
Put $\varrho=\mu (2\Theta-v)$, one has 
\begin{align*}
	\lambda_t\leq \varrho t.
\end{align*}
In particular, by the definition of $\varrho$ and $\theta$ one obtains that, for any $t\leq\aleph$, 
\begin{align}
	\lambda_t\leq \theta \aleph. 
	\label{queu-bound-loc}
\end{align}
Take now $T\geq \aleph$ and suppose that the bound in $(\ref{queu-bound-loc})$ holds for any $t\leq T$. Our goal is to show that it also holds at time $T+1$ and thus, by induction, at any time $t\geq T$. Notice that, if $\lambda_{T+1}\leq \lambda_{T-\aleph}$, then $\lambda_{T+1}\leq \theta \aleph$ immediately. Therefore, we suppose without loss of generality that $\lambda_{T+1}> \lambda_{T-\aleph}$. Summing the inequality in $(\ref{Delt-bound})$ over $t\in\{T-\aleph,\ldots,T \}$ gives 
\begin{equation}
	\begin{split}    
		\sum_{t=T-\aleph}^{T} \frac{\Delta(\lambda_t)}{\mu}&\leq  (\aleph+1) \frac{B}{\mu}+(\aleph+1)\frac{\alpha \Theta^2}{4}-\eta \sum_{t=T-\aleph}^{T} \lambda_t\\
		&+\frac{1}{2\alpha}\sum_{t=T-\aleph}^{T} \bigg(\|Q^*-P^{t-1}\|^2-\|Q^*-P^{t}\|^2\bigg)\\
		&+\sum_{t=T-\aleph}^{T} \Theta. 
	\end{split}
\end{equation} 
Therefore
\begin{equation}
	\begin{split}    
		\frac{1}{\mu}(\lambda^2_{t+1}-\lambda^2_{T-\aleph})&\leq  (\aleph+1) \frac{B}{\mu}+(\aleph+1)\frac{\alpha \Theta^2}{4}-\eta \sum_{t=T-\aleph}^{T} \lambda_t\\
		&+\frac{1}{2\alpha} +(\aleph+1)\Theta. 
	\end{split}
\end{equation} 
Recall that by assumption  $\lambda_{T+1}> \lambda_{T-\aleph}$. Therefore 
\begin{equation}
	\begin{split}    
		\eta \sum_{t=T-\aleph}^{T} \lambda_t&\leq  (\aleph+1) \frac{B}{\mu}+(\aleph+1)\frac{\alpha \Theta^2}{4}+\frac{1}{2\alpha} +(\aleph+1)\Theta.  
		\label{loc-sum-qu-ineq}
	\end{split}
\end{equation}
In order to prove that the inequality in $(\ref{queu-bound-loc})$ holds at time $t=T+1$, we proceed by contradiction. Indeed, suppose that 
\begin{align}
	\lambda_{T+1}> \theta \aleph. 
	\label{contr-hyp}
\end{align}
From $(\ref{lamb-succ-bound})$ one obtains, for all $t\in\{0,1,\ldots,T+1\}$,
\begin{align*}
	\lambda_{T+1}-\lambda_t\leq (T+1-t)\varrho. 
\end{align*}
Thus, the two last inequalities give   
\begin{align*}
	\lambda_t\geq \theta\aleph-(T+1-t)\varrho. 
\end{align*}
Replacing $\lambda_t$ in (\ref{loc-sum-qu-ineq}) by the last upper-bound leads to 
\begin{equation}
	\begin{split}    
		\eta \sum_{t=T-\aleph}^{T} (\theta\aleph-(T+1-t)\varrho)&\leq  (\aleph+1) \frac{B}{\mu}+(\aleph+1)\frac{\alpha \Theta^2}{4}\\
		&+\frac{1}{2\alpha} +(\aleph+1)\Theta.  
		\label{loc-sum-qu-ineq}
	\end{split}
\end{equation}
Developing the sum gives 
\begin{equation} 
	\begin{split}
		\eta(\aleph+1)(\theta\aleph-\varrho\frac{(\aleph+2)}{2})&\leq  (\aleph+1) \frac{B}{\mu}+(\aleph+1)\frac{\alpha \Theta^2}{4}\\
		&+\frac{1}{2\alpha} +(\aleph+1)\Theta.  
		\label{loc-sum-qu-ineq}
	\end{split}
\end{equation}
Thus
\begin{equation} 
	\begin{split}
		\theta&\leq \frac{1}{\eta\aleph}\bigg(  \frac{B}{\mu}+  \frac{\alpha \Theta^2}{4}+\frac{1}{(\aleph+1)}  \frac{1}{2\alpha} +  \Theta\bigg)+ \frac{\varrho(\aleph+2)}{2\aleph}.  
		\label{loc-sum-qu-ineq}
	\end{split}
\end{equation}
However, this contradicts the definition of $\theta$ in $(\ref{thet-def})$. Therefore, the inequality in $(\ref{contr-hyp})$ cannot be true! We thus deduce that 
\begin{align}
	\lambda_{T+1}\leq \theta \aleph,
\end{align}
and, by induction that, for all $t\geq T$, $\lambda_t\leq \theta \aleph$.
\carre

\subsection{Regret against $K$ benchmark} 
\label{K-bench-proof}
Using Lemma \ref{Delta-bound} with $P=P_*^K$ one obtains  that for every slot $t \in\{1, 2,\ldots\}$, the solution $P^t$ given by Algorithm $\ref{sad-algo}$ satisfies 
\begin{align*}  
	\frac{\Delta(\lambda_t)}{\mu}&\leq \frac{B}{\mu}+ \mathbb{E}_{P_*^K}[C_R(A^t)]-\mathbb{E}_{P^{t-1}}[C_R(A^{t-1})]\\
	&+\lambda_t \big(\mathbb{E}_{P_*^K}[C(A^t,b^{t-1})\big]-v \big)\\
	&+\frac{1}{2\alpha} \bigg(\|P_*^K-P^{t-1}\|^2-\|P_*^K-P^{t}\|^2\bigg)+\frac{\alpha \Theta^2}{4}.
\end{align*}
Summing the last inequality over the slots $\{t,\ldots,t+K-1 \}$ gives
\begin{align*}  
	&\sum_{\tau=0}^{K-1} \frac{\Delta(\lambda_{t+\tau})}{\mu}\\
	&\leq \frac{BK}{\mu}+K\frac{\alpha \Theta^2}{4}\\
	&\quad+\sum_{\tau=0}^{K-1}\big( \mathbb{E}_{P_*^K}[C_R(A)]-\mathbb{E}_{P^{t+\tau-1}}[C_R(A)]\big)\\
	&\quad+\sum_{\tau=0}^{K-1} \lambda_{t+\tau} \big(\mathbb{E}_{P_*^K}[C(A,b^{t+\tau-1})\big]-v \big)\\
	&\quad+\frac{1}{2\alpha} \sum_{\tau=0}^{K-1} \bigg(\|P_*^K-P^{t+\tau-1}\|^2-\|P_*^K-P^{t+\tau}\|^2\bigg).
\end{align*}
Moreover, note that
\begin{align*}
	&\sum_{\tau=0}^{K-1} \lambda_{t+\tau} \bigg(\mathbb{E}_{P^K_*}\big[C(A,b^{t+\tau-1})\big]-v \bigg)\\
	&= \sum_{\tau=0}^{K-1} \lambda_{t+\tau}\bigg[ \big(\mathbb{E}_{P^K_*}\big[C(A,b^{t+\tau-1})\big]-v \big)^+\\
	&\qquad\qquad\qquad+\big(\mathbb{E}_{P^K_*}\big[C(A,b^{t+\tau-1})\big]-v \big)^-\bigg].
\end{align*}
The equation in $(\ref{lambd-update})$ gives, for any $\tau\geq 1$,   
\begin{align*}
	\lambda_{t+\tau}\geq \lambda_{t}+ \mu \sum_{s=0}^{\tau-1} \big(\mathbb{E}_{P_{t+s}}\big[C(A,b^{t+s-1})\big]-v \big).
\end{align*}
Moreover, from the assumptions above one gets 
\begin{align*}
	\lambda_{t+\tau}\leq \lambda_{t}+\mu \sum_{s=0}^{\tau-1}\big(2\Theta-v \big).
\end{align*}
Thus,  
\begin{align*}
	&\sum_{\tau=0}^{K-1} \lambda_{t+\tau} \big(\mathbb{E}_{P^K_*}\big[C(A,b^{t+\tau-1})\big]-v \big)\\
	&\leq \sum_{\tau=0}^{K-1} \bigg(\lambda_{t}+\mu \sum_{s=0}^{\tau-1}\big(2\Theta-v \big)\bigg)  \big(\mathbb{E}_{P^K_*}\big[C(A,b^{t+\tau})\big]-v \big)^+\\
	&\qquad+ \bigg(\lambda_{t}+ \mu \sum_{s=0}^{\tau-1} \big(\mathbb{E}_{P_{t+s}}\big[C(A,b^{t+s-1})\big]-v \big) \bigg)\\
	&\qquad\qquad\times\big(\mathbb{E}_{P^K_*}\big[C(A,b^{t+\tau})\big]-v \big)^-\\ 
	&\leq \lambda_{t}\sum_{\tau=1}^{K-1}  \big(\mathbb{E}_{P^K_*}\big[C(A,b^{t+\tau})\big]-v \big) +\frac{1}{2} K(K-1)\mu \big(2 \Theta-v \big)^2\\
	&\leq \frac{1}{2} K(K-1)\mu \big(4 \Theta^2+v^2 \big), 
\end{align*} 
where the last inequality follows by the fact that the constraint is satisfied over any finite $K$ slots. Therefore, 
\begin{align*}  
	\\& \sum_{\tau=0}^{K-1} \frac{\Delta(\lambda_{t+\tau})}{\mu}\\
	&\quad\leq \sum_{\tau=0}^{K-1}\big( \mathbb{E}_{P_*^K}[C_R(A)]-\mathbb{E}_{P^{t+\tau-1}}[C_R(A)]\big)\\
	&\qquad+\frac{1}{2}K^2\mu \big(4 \Theta^2+v^2 \big)\\
	&\qquad+\frac{1}{2\alpha} \sum_{\tau=0}^{K-1} \bigg(\|P_*^K-P^{t+\tau-1}\|^2-\|P_*^K-P^{t+\tau}\|^2\bigg)\\
	&\qquad+\sum_{\tau=0}^{K-1}\frac{\alpha \Theta^2}{4}.
\end{align*} 
Summing the last inequality for $t=\{1,2,\ldots, T-K\}$ gives 
\begin{align*}  
	&\sum_{\tau=0}^{K-1} \sum_{t=1}^{T-K} \frac{\Delta(\lambda_{t+\tau})}{\mu}\\
	&\leq \sum_{\tau=0}^{K-1}\sum_{t=1}^{T-K}\big( \mathbb{E}_{P_*^K}[C_R(A)]-\mathbb{E}_{P^{t+\tau-1}}[C_R(A)]\big)\\
	&+\sum_{t=1}^{T-K}\frac{1}{2} K^2\mu \big(4 \Theta^2+v^2 \big)\\
	&+\frac{1}{2\alpha} \sum_{\tau=0}^{K-1}\sum_{t=1}^{T-K} \bigg(\|P_*^K-P^{t+\tau-1}\|^2-\|P_*^K-P^{t+\tau}\|^2\bigg)\\
	&+\sum_{t=1}^{T-K}K\frac{\alpha \Theta^2}{4}. 
\end{align*} 
By the update equation in $(\ref{lambd-update})$ together with Assumption \ref{ass}, one gets   
\begin{align*}
	\sum_{\tau=0}^{K-1} \sum_{t=1}^{T-K} \frac{\Delta(\lambda_{t+\tau})}{\mu}&=\frac{1}{\mu} \sum_{\tau=0}^{K-1} \sum_{t=1}^{T-K} (\lambda^2_{t+\tau+1}-\lambda^2_{t+\tau})\\
	&=\frac{1}{\mu} \sum_{\tau=0}^{K-1}  (\lambda^2_{T-K+\tau+1}-\lambda^2_{1+\tau})\\
	&\geq -\frac{1}{\mu} \sum_{\tau=0}^{K-1}  (\lambda^2_{\tau+1})\\ 
	&\geq -\frac{1}{\mu} \sum_{\tau=0}^{K-1}  (\tau \mu(2\Theta+v))^2\\ 
	&\geq - \mu(2\Theta+v)^2 \frac{(K-1)K(2K-1)}{6}. 
\end{align*}
Moreover, simple calculations allows to obtain 
\begin{align*}
	\frac{1}{2} \sum_{\tau=0}^{K-1}\sum_{t=1}^{T-K} \big(\|P_*^K-P^{t+\tau-1}\|^2-\|P_*^K-P^{t+\tau}\|^2\big)\leq K.
\end{align*}
Then, 
\begin{align*}  
	&\sum_{\tau=0}^{K-1}\sum_{t=1}^{T-K}\big(\mathbb{E}_{P^{t+\tau-1}}[C_R(A)] -\mathbb{E}_{P_*^K}[C_R(A)]\big)\\
	&\leq\mu(2\Theta+v)^2 \frac{(K-1)K(2K-1)}{6} \\
	&+\sum_{t=1}^{T-K}\frac{1}{2} K^2\mu \big(4 \Theta^2+v^2 \big)+\frac{1}{\alpha} K+\sum_{t=1}^{T-K}K \frac{\alpha \Theta^2}{4}.
\end{align*} 
Note that, for any $1\leq K\leq T$, it is easy to verify that  
\begin{align*}
	&K\sum_{t=1}^T\big(\mathbb{E}_{P^t}[C_R(A)]-\mathbb{E}_{P^K_*}[C_R(A)]\big)\\
	&=\sum_{\tau=0}^{K-1}\sum_{t=1}^{T-K}\big(\mathbb{E}_{P^{t+\tau}}[C_R(A)]-\mathbb{E}_{P^K_*}[C_R(A)]\big)\\
	&\quad+\sum_{\tau=0}^{K-1}\sum_{t=1}^{\tau}\big(\mathbb{E}_{P^t}[C_R(A)]-\mathbb{E}_{P^K_*}[C_R(A)]\big)\\
	&\quad+\sum_{\tau=0}^{K-1}\sum_{t=T-K+\tau+1}^{T}\big(\mathbb{E}_{P^t}[C_R(A)]-\mathbb{E}_{P^K_*}[C_R(A)]\big). 
\end{align*}
Therefore,
\begin{align*}
	&K\sum_{t=1}^T\big(\mathbb{E}_{P^t}[C_R(A)]-\mathbb{E}_{P^K_*}[C_R(A)]\big)\\
	&\leq \mu(2\Theta+v)^2 \frac{(K-1)K(2K-1)}{6} \\
	&+\sum_{t=1}^{T-K}\frac{1}{2} K^2\mu \big(4 \Theta^2+v^2 \big)+\frac{1}{\alpha} K+\sum_{t=1}^{T-K}K \frac{\alpha \Theta^2}{4}\\
	&+ \frac{K(K-1)}{2}\Theta+ K(K-1)\Theta-\frac{K(K-1)}{2}\Theta. 
\end{align*}
Thus,
\begin{align*}
	&\sum_{t=1}^T\big(\mathbb{E}_{P^t}[C_R(A)]-\mathbb{E}_{P^K_*}[C_R(A)]\big)\\
	&\quad\leq \mu(2\Theta+v)^2 \frac{(K-1)(2K-1)}{6} \\
	&\qquad+\sum_{t=1}^{T-K}\frac{1}{2} K\mu \big(4 \Theta^2+v^2 \big)+\frac{1}{\alpha} \\
	&\qquad+\sum_{t=1}^{T-K} \frac{\alpha \Theta^2}{4}+  K(K-1)\Theta,
\end{align*}
which concludes the proof.

\subsection{Proof of Corollary \ref{coro}} 
\label{coro-proof}
Define the random variables 
\begin{align*}
	X^t= C_R(A^t)-\mathbb{E}_{P^t}[C_R(A^t)]. 
\end{align*}
Then, by Assumption \ref{ass}, $|X^t|\leq \Theta$. Therefore, $\{X_t\}_{t\geq 1}$ is a sequence of bounded martingales differences, and $Y_T=\sum_{t=1}^TX_t$ is a martingale with respect to the filtration $\mathcal{F}_T=\sigma(A^t,t\leq T)$. Thus, Hoeffding-Azuma's inequality (see e.g. \cite{b23}) gives us that, for any $0<\delta<1$, 
\begin{align*}
	\mathbb{P}\bigg(Y_T\leq \sqrt{2 \log (\delta^{-1}) T \Theta^2}\bigg)\geq 1-\delta.
\end{align*}  
Thus, with probability at least $1-\delta$, one has 
\begin{align*}
	\sum_{t=1}^T C_R(A^t) \leq \sum_{t=1}^T \mathbb{E}_{P^t} [C_R(A^t)]+\sqrt{2 \log (\delta^{-1}) T \Theta^2}. 
\end{align*} 
Now, it is easy to see that $\sum_{t=1}^T\mathbb{E}_{P_K^*}[C_R(A^t)]\leq \sum_{t=1}^T C_R(A^*)$, where $P^K_*$ is the solution to the problem in $(\ref{p-star-stat-K})$ (since $\delta_{A^*}\in\mathcal{P}(\mathcal{R})$). Therefore, 
\begin{equation}
	\begin{split}
		\tilde{R}^{K}_T &\leq R^{K}_T+\sqrt{2 \log (\delta^{-1}) T \Theta^2}. 
	\end{split}
\end{equation}
Finally, using Theorem \ref{Theo-regret-K} gives (\ref{reg-non-rand-bound})
\carre

\vspace{12pt}
\color{red}

\end{document}